\input amstex
\documentstyle{amsppt}

\magnification=\magstephalf

\hsize=30pc
\vsize=42pc
\hoffset=.5truein
\voffset=.75truein

\document
\refstyle{A}

\topmatter

\title
$X$-Inner Automorphisms of Semi-Commutative Quantum Algebras
\endtitle

\rightheadtext{$X$-Inner Automorphisms of Semi-Commutative Quantum Algebras}

\leftheadtext{Jeffrey Bergen and Mark C. Wilson}

\author
Jeffrey Bergen and Mark C. Wilson
\endauthor

\address
Department of Mathematics, DePaul University,
Chicago, Illinois 60614
\endaddress

\email
jbergen\@condor.depaul.edu
\endemail

\address
Department of Mathematics, University of Auckland,
Private Bag 92019, Auckland, New Zealand
\endaddress

\email
wilson\@math.auckland.ac.nz
\endemail

\thanks
Much of this work was done when the first author was a visitor at the
University of Auckland and he would like to thank the department for
its hospitality.
The second author was supported by a NZST Postdoctoral Fellowship.
\endthanks

\subjclass
17B37, 16W20, 16S36, 16S30
\endsubjclass

\abstract
Many important quantum algebras such as quantum symplectic space, quantum Euclidean space,
quantum matrices, $q$-analogs of the Heisenberg algebra and the quantum Weyl algebra are
semi-commutative.  In addition, enveloping algebras $U(L_+)$ of even Lie color algebras are also
semi-commutative.  In this paper, we generalize work of Montgomery and examine the $X$-inner
automorphisms of such algebras.

The theorems and examples in our paper show that for algebras $R$ of this type, the non-identity
$X$-inner automorphisms of $R$ tend to have infinite order.  Thus if 
$G$ is a finite group of automorphisms of $R$, then the action of $G$ will be $X$-outer and this
immediately gives  useful information about crossed products $R*_tG$.
\endabstract

\endtopmatter

\head
\S 1. Introduction and Terminology
\endhead 

McConnell \cite{Mc} and Musson \cite{Mu} have shown that many important quantum algebras are
semi-commutative.  In this paper we examine the $X$-inner automorphisms of such algebras.
In particular, knowing the $X$-inner automorphisms of these algebras is helpful when studying
crossed products.

In several previous papers, the $X$-inner automorphisms of some naturally
occuring prime rings are computed.  In \cite{M}, Montgomery computes the
$X$-inner automorphisms of filtered algebras $R$ whose associated
graded ring $gr(R)$ is a commutative domain.  Montgomery then uses this result
to classify the $X$-inner automorphisms of enveloping algebras $U(L)$ of Lie
algebras $L$, as well as the $R$-stabilizing $X$-inner automorphisms of
skew polynomial rings $R[x;d]$, where $d$ is a derivation of $R$. 

The results in \cite{M} have been generalized in several directions.  In 
\cite{OP}, Osterburg and Passman compute the $R$-stabilizing $X$-inner
automorphisms of smash products $R\#U(L)$.  In \cite{LM}, Leroy and
Matczuk compute the $R$-stabilizing $X$-inner automorphisms of skew polynomial 
rings $R[x;\tau, \delta]$, where $\tau$ is an automorphism of $R$ and $\delta$ 
is a $\tau$-derivation of $R$.  

In \cite{Mc}, McConnell defines a {\it semi-commutative} algebra
as a filtered algebra $R$ whose associated graded ring $gr(R)$
is the coordinate ring of quantum space $\Cal O_q(K^n)$. 
We will be somewhat informal and simple refer to $\Cal O_q(K^n)$ as {\it quantum space}.
An example of a semi-commutative algebra is the quantum Weyl algebra $K[x,y \mid xy-qyx=1]$ whose
associated  graded ring is the quantum plane.

McConnell and Musson show that many
important quantum algebras such as quantum symplectic space, quantum Euclidean space and
quantum matrices can be filtered in such a way that they are indeed semi-commutative.  
In addition, enveloping algebras $U(L_+)$ of even Lie color algebras are also semi-commutative.
In light of these examples,  it is natural to try to extend Montgomery's work on $X$-inner
automorphisms to semi-commutative algebras.

In Section 2, we begin by computing the $X$-inner automorphisms of a
generalized version of quantum space. This enables us to obtain information
about the $X$-inner automorphisms of semi-commutative algebras.
Next, we compute the $X$-inner automorphisms of enveloping algebras $U(L_+)$ of 
even Lie color algebras $L_+$.   
This allows us to analyze the structure of the homogeneous
semi-invariants under the adjoint action of $L_+$ in the Martindale quotient ring
of $U(L_+)$.  In particular, we show that every $X$-inner automorphism of
$U(L_+)$ is induced by a homogeneous semi-invariant.   
In addition, we show that the group of $X$-inner automorphisms of $U(L_+)$ is abelian and the group of homogeneous semi-invariants is nilpotent of class $2$.
We conclude this section by applying our results to the 
study of crossed products $U(L_+)*_tG$, where $G$ is a finite group of automorphisms of $U(L_+)$.

In Section 3, we look at the $(R,R)$-bimodules of $R[x;\tau, \delta]$ in the
special case where $\delta$ is a $q$-skew derivation.  This enables us to 
sharpen a result in \cite{LM} on the $X$-inner automorphisms of 
$R[x;\tau,\delta]$.

Finally, in Section 4 we apply the results of the previous two sections to compute the
$X$-inner automorphisms of various quantum algebras.  
We first consider the quantum Weyl algebra $K[x,y \mid xy-qyx=1]$ and then the enveloping algebra
of a $3$-dimensional even Lie color algebra.  
Next, we look at two algebras studied in \cite{KS} by Kirkman and Small which are $q$-analogs of
the enveloping algebra of the Heisenberg algebra.  Although the generators and relations of these
two algebras are almost identical, we show that for one of them their group of $X$-inner
automorphisms is infinite cyclic, whereas the group of $X$-inner automorphisms
of the other is free abelian on two generators.  We then compute the $X$-inner automorphisms of an algebra which, although quite different
from $q$-analogs of enveloping algebras and Weyl algebras, is indeed semi-commutative and therefore can be studied using the techniques developed in the previous sections. We conclude this paper by observing that all of the
non-identity $X$-inner automorphisms in our examples have infinite order and this has immediate applications to the structure of crossed products.

At this point, we introduce some of the terminology that we will use throughout this
paper. If $R$ is a prime ring, we let $Q(R)$ denote its symmetric Martindale
quotient ring and let $N(R)$ denote those $w \in Q(R)$ such that $wR = Rw$. 
An automorphism $\sigma$ of $R$ is said to be {\it $X$-inner} if it induced by some element in $Q(R)$.  In particular, this means that there is some $w \in N(R)$ such that $r^\sigma w  = wr$, for all $r \in R$.  
An important and well known fact is that an automorphism $\sigma$ is
$X$-inner if and only if there exist nonzero $a,b \in R$ such that
$ar^\sigma b^\sigma = bra$, for all $r \in R$. 

If $\tau$ is an automorphism of $R$, then an additive map $\delta: R \longrightarrow R$
is said to be a { \it skew derivation} of $R$ if
$\delta(rs) = \delta(r) s + \tau(r) \delta(s)$, for all $r,s \in R$.  In addition, we say that
$\delta$ is a {\it $q$-skew derivation} if it is also that case that
$\delta \tau = q \tau \delta$, where $R$ is an algebra over a field $K$ and
$q$ is a nonzero element of $K$. A skew derivation $\delta$ is said to be
$X$-inner if there is some $c \in Q(R)$ such that $\delta(r) = cr - \tau(r) c$, for all
$r \in R$.

\head
\S 2. Quantum Algebras and Lie Color Algebras
\endhead

Let $X = \{x_i \mid i \in I \}$ be a set of variables
indexed by the set $I$.  
If $K$ is a field, we define the algebra $K_q[X]$ to be the $K$-algebra generated
by the elements of the set $X$ subject to the conditions that
$x_ix_j = q_{ij}x_jx_i$, for all $i,j \in I$, where the $q_{ij}$ are nonzero
elements of $K$. 
Clearly we must have that $q_{ii} = 1$ and  $q_{ij} = {q_{ji}}^{-1}$, 
however we impose no additional conditions on the  $q_{ij}$.
Note that we are not necessarily assuming that the set $I$ is finite; however
if $I$ has only two elements then $K_q[X]$ is the quantum plane,
whereas if $I$ has $n$ elements then $K_q[X]$ is a generalized version of quantum
$n$-space.  

When studying $K_q[X]$ it will also be useful to consider the larger
$K$-algebra $K_q[X,X^{-1}]$ which is generated by the elements of $X$ as well as their
inverses subject to the condition that $x_ix_j = q_{ij}x_jx_i$.  
These algebras are twisted group rings of free abelian groups and have been studied 
in \cite{McP} by McConnell and Pettit and in \cite{OsbP} by Osborn and Passman.

If $\Delta$ is a monomial in $K_q[X,X^{-1}]$ then, for every $x \in X$,
there is a nonzero $\beta_x \in K$ such that $\Delta x \Delta^{-1} = \beta_x x$.
Thus conjugation by $\Delta$ induces an automorphism of $K_q[X]$.
Note that there exist monomials $\Delta_1 , \Delta_2 \in  K_q[X]$ and a nonzero
$\gamma \in K$ such that $\gamma \Delta = \Delta_1^{-1}\Delta_2$.  If we let
$\sigma$ denote the automorphism induced by $\Delta$, then 
$\Delta_1^{-1}\Delta_2$ also induces $\sigma$.
Hence
$$
s^\sigma = (\Delta_1^{-1}\Delta_2)s(\Delta_1^{-1}\Delta_2)^{-1},
$$
for all $s \in K_q[X]$.
In the above equation if we replace $s$ by $s\Delta_2$ and multiply on the left by
$\Delta_1$, we obtain
$$
\Delta_1s^\sigma\Delta_2^\sigma = \Delta_2s\Delta_1,
$$
and so, $\sigma$ is an $X$-inner automorphism of $K_q[X]$.
In our first result, we will show that this is how all $X$-inner automorphisms of
$K_q[X]$ arise.

For any monomial $\Delta$ in $K_q[X]$ and $x \in X$, we will let $\pi_\Delta(x)$ denote
the element of $K$ such that 
$$
\Delta x \Delta^{-1} = \pi_\Delta(x) x.
$$
Thus we can consider $\pi_\Delta$ as a map from $X$ to $K$.
We say that a nonzero element $w \in K_q[X]$ is {\it monotone} if
$\pi_\Delta$ is the same for every monomial $\Delta$ in the support of $w$.
Note that saying $w$ is monotone is equivalent to saying that every monomial in
the support of $w$ induces the same $X$-inner automorphism of $K_q[X]$.

\proclaim{Lemma 2.1}
{\rm (a)}  Every nonzero ideal of $K_q[X]$ contains a monotone element.
\roster
\item"(b)"  Every monotone element of $K_q[X]$ is of the form
$\Delta f$, where $\Delta$ is a monomial in $K_q[X]$ and $f$ belongs to
the center of $K_q[X,X^{-1}]$.
\item"(c)"  An element of $K_q[X]$ is normalizing if and only if it is  monotone.  Thus every nonzero ideal of $K_q[X]$ contains a normalizing element. 
\endroster
\endproclaim

\demo{Proof}
(a)  If $J$ is a nonzero ideal of $K_q[X]$, let $w$ be a nonzero element of 
$J$ with a minimal number of monomials in its support.  We claim that
$w$ is indeed monotone.

To this end, if $w = \sum_i \gamma_i\Delta_i$ and $x \in X$, then we consider
$$
v = \pi_{\Delta_1}(x)xw -wx = \sum_i\gamma_i(\pi_{\Delta_1}(x)x\Delta_i -\Delta_ix) =
\sum_i \gamma_i(\pi_{\Delta_1}(x) - \pi_{\Delta_i}(x)) x\Delta_i.
$$
Since $x\Delta_1$ is not in the support of $v$, it follows that $v$ has fewer elements
in its support then does $w$.  Furthermore, since $v \in J$ we see that $v=0$ and thus
$\pi_{\Delta_1}(x) = \pi_{\Delta_i}(x)$, for all $x \in X$.  Hence 
$\pi_{\Delta_1} = \pi_{\Delta_i}$, for all $i$, and $w$ is indeed monotone.

(b)  If $w \in K_q[X]$ is monotone and if $\Delta$ is a monomial
in the support of $w$, then every monomial in $w$ induces the same $X$-inner
automorphism as $\Delta$.  Thus $wxw^{-1} = \pi_\Delta(x)x$, for all $x \in X$.
Hence $w$ induces the same $X$-inner automorphism as $\Delta$ and this implies that
$\Delta^{-1}w$ must be central in $K_q[X,X^{-1}]$.  Therefore if we let
$f = \Delta^{-1}w$, then we have $w = \Delta f$, as desired.

(c)  It is clear that every monotone element is normalizing.  Now, if $b \in K_q[X]$ is normalizing, let $\tau$ denote the $X$-inner
automorphism of $K_q[X]$ induced by $b$.  Thus $s^\tau b = bs$, for all $s \in K_q[X]$.
Let $N$ and $M$ denote, respectively, the largest and smallest degrees of monomials in
the support of $b$.  If $x \in X$, then clearly $N+1$ and $M+1$ are, respectively, the
largest and smallest degrees of monomials in the support of $bx$.
However, $x^\tau b = bx$, thus $N+1$ and $M+1$ must also be, respectively, the
largest and smallest degrees of monomials in the support of $x^\tau b$.
Thus $x^\tau$ must be a linear combination of elements of $X$.
However, if $x^\tau$ includes some $x^\prime \in X$ other than $x$ in its support,
then the degree of $x^\prime$ in $x^\tau b$ will exceed its degree in $bx$, a
contradiction.  Thus $x^\tau = \gamma_xx$, for some nonzero $\gamma_x \in K$.

Finally, if $b = \sum_i \gamma_i \Delta_i$ then
$$
x^\tau b = \gamma_x x \sum_i \gamma_i \Delta_i = \sum_i \gamma_i \gamma_x x \Delta_i
$$ 
and
$$
bx = (\sum_i \gamma_i \Delta_i)x = \sum_i \gamma_i \Delta_i x =
\sum_i\gamma_i\pi_{\Delta_i}(x)x\Delta_i.
$$
Thus $\pi_{\Delta_i}(x) = \gamma_x$, for all $x \in X$.  As a result
$\pi_{\Delta_i}$ is the same for every $\Delta_i$ in the support of $b$ and so,
$b$ is monotone. 
It now follows, from part (a), that every nonzero ideal contains a normalizing element.
\qed
\enddemo

We can now describe the $X$-inner automorphisms of quantum space.

\proclaim{Proposition 2.2}
Every $X$-inner automorphism $\sigma$ of $K_q[X]$ is induced by conjugation by some
monomial $\Delta$ in $K_q[X,X^{-1}]$.  In particular, if
$\Delta = x_{i_1}^{j_1}x_{i_2}^{j_2}\cdots x_{i_n}^{j_n}$ then
${x_k}^\sigma = \gamma x_k$ where 
$\gamma = q_{{i_1}k}^{j_1}q_{{i_2}k}^{j_2}\cdots q_{{i_n}k}^{j_n}$, for all $k \in I$. 
\endproclaim

\demo{Proof}
Let $q \in Q(K_q[X])$ such that $q$ induces $\sigma$, thus
$$
s^\sigma q = qs,
\tag 1
$$
for all $s \in K_q[X]$.
Let $J \not=0$ be an ideal of $K_q[X]$ such that $Jq, qJ \subseteq K_q[X]$.
By Lemma 2.1(a), there exists a monotone element $a \in J$.
If $b = aq$, then $b$ is also a normalizing element of $K_q[X]$, thus by
Lemma 2.1(c), $b$ is also monotone.

By multiplying equation (1) on the left and right by $a$ we see that
$as^\sigma qa = aqsa$.  Since $b = aq$, it follows that
$b^\sigma = qa$ and we obtain 
$$
as^\sigma b^\sigma = bsa.
\tag 2
$$
By Lemma 2.1(b), $a = \Delta_1f$ and $b = \Delta_2g$, where $\Delta_1, \Delta_2$ are 
monomials in $K_q[X]$ and $f,g$ are in the center of $K_q[X,X^{-1}]$.
Since $g^\sigma = g$, equation (2) becomes
$\Delta_1fs^\sigma \Delta_2^\sigma g = \Delta_2gs\Delta_1f$, which reduces to
$\Delta_1s^\sigma \Delta_2^\sigma = \Delta_2 s \Delta_1$.
Thus
$$
s^\sigma = \Delta_1^{-1}\Delta_2s\Delta_1(\Delta_2^\sigma)^{-1}.
\tag 3
$$

Letting $s=1$ in equation (3) yields $\Delta_1\Delta_2^\sigma = \Delta_2 \Delta_1$,
which implies that $(\Delta_2^\sigma)^{-1} = \Delta_1^{-1}\Delta_2^{-1}\Delta_1$.
Plugging this fact into equation (3) results in
$$
s^\sigma = \Delta_1^{-1}\Delta_2s\Delta_1\Delta_1^{-1}\Delta_2^{-1}\Delta_1 =
(\Delta_1^{-1}\Delta_2)s(\Delta_1^{-1}\Delta_2)^{-1}.
$$
Thus $\sigma$ is induced by the monomial $\Delta = \Delta_1^{-1}\Delta_2 
\in K_q[X,X^{-1}]$.

Finally, if 
$\Delta = x_{i_1}^{j_1}x_{i_2}^{j_2}\cdots x_{i_n}^{j_n}$ then
for any $x_k \in X$, we have
$$
{x_k}^\sigma = \Delta x_k \Delta^{-1} = 
x_{i_1}^{j_1}x_{i_2}^{j_2}\cdots x_{i_n}^{j_n}x_k
x_{i_n}^{-j_n}\cdots x_{i_2}^{-j_2}x_{i_1}^{-j_1} =
q_{{i_1}k}^{j_1}q_{{i_2}k}^{j_2}\cdots q_{{i_n}k}^{j_n}x_k.
$$
\qed
\enddemo

In \cite{Mc}, McConnell shows that many quantum algebras $R$ can be filtered in 
such a way that their associated graded algebra $gr(R)$ is isomorphic to $K_q[X]$.
We will use Proposition 2.2 to obtain information about the $X$-inner automorphisms of
such algebras.

If $R$ is a filtered $K$-algebra, then $R = \bigcup_{n \ge0} R_n$, where
$K \subseteq R_0$.  If $r \in R$, then $d(r)$ will denote the filtration of
$r$, thus $d(r) = m$ if $r \in R_m$ and $r \notin R_{m-1}$.
The associated graded algebra $gr(R)$ is $\oplus_{n \ge0} R_n/R_{n-1}$.
If $d(r) = m$, then we let ${\overline r}$ denote the element $r + R_{m-1}$ in $gr(R)$. 
We can now prove our first main result, which is an analog of Montgomery's result
\cite{M} on the $X$-inner automorphisms of filtered algebras whose associated graded
ring is a commutative domain.

\proclaim{Theorem 2.3}
Let $R$ be a filtered $K$-algebra such that the associated graded algebra $gr(R)$ is
isomorphic to $K_q[X]$.  If $\sigma$ is an $X$-inner automorphism of $R$, then
$\sigma$ preserves the filtration of $R$.  Furthermore, if $x \in R$
such that 
${\overline x} \in X$ then $x^\sigma = \alpha_x x +$ terms of lower degree,
where $\alpha_x$ is a nonzero element of $K$.
\endproclaim

\demo{Proof}
Since $\sigma$ is an $X$-inner automorphism of $R$, there exists nonzero $a,b \in R$
such that $as^\sigma b^\sigma = bsa$, for all $s \in R$.
Letting $s=1$, we have $ab^\sigma = ba$.  Since $gr(R)$ has an additive degree function,
we have $d(a) + d(b) = d(ab^\sigma) = d(ba) = d(b) + d(a)$.  Thus $d(b^\sigma) = d(b)$.

It now follows that, for any $s \in R$, 
$d(a) + d(s^\sigma) + d(b) = d(a) + d(s^\sigma) + d(b^\sigma) = d(as^\sigma b^\sigma) =
d(bsa) = d(b) + d(s) + d(a)$.  Thus $d(s^\sigma) = d(s)$ and so, $\sigma$ preserves the
filtration of $R$.

Since $\sigma$ does preserve the filtration, $\sigma$ induces an automorphism
${\overline \sigma}$ of $gr(R)$ defined as
${\overline \sigma}(r + R_{m-1}) = \sigma(r) + R_{m-1}$, for all $r \in R_m$.
The equation $as^\sigma b^\sigma = bsa$ becomes
${\overline a} {\overline {s^\sigma}} {\overline {b^\sigma}} = {\overline b} {\overline s} {\overline a}$,
thus ${\overline \sigma}$ is $X$-inner on $gr(R)$.
By Proposition 2.2, if $x \in R$ such that ${\overline x} \in X$ then
${\overline \sigma}({\overline x}) = \alpha_x {\overline x}$, for some nonzero $\alpha_x \in K$.
Thus $\sigma(x) - \alpha_x x$ has lower degree than $x$, hence
$\sigma(x) = \alpha_x x +$ terms of lower degree.
\qed
\enddemo

We now briefly look at the type of algebra to which Theorem 2.3 can be applied.
\vskip.1in

\noindent {\bf Example 2.4.}
Let $R$ be the $K$-algebra generated by the set $\{x,y,z\}$ subject to the relation
$xy - qyx = \alpha yz + \beta y^2 + \gamma y$, where $q, \alpha, \beta, \gamma \in K$
and $q$ is nonzero.
We can filter $R$ by letting $d(x) = 2, d(y) = 1, d(z) = 1$.
Therefore $gr(R)$ is generated by ${\overline x}, {\overline y}, 
{\overline z}$ subject to the
sole relation ${\overline x} {\overline y} = q {\overline y} {\overline x}$.  Thus
$gr(R)$ is quantum $3$-space.  Therefore Theorem 2.3 implies that any $X$-inner
automorphism of $R$ must send $x$ to a scalar multiple of $x$ plus terms of lower degree.
For an example of such an automorphism, note that $yR = Ry$ and so,  $y$ induces an
$X$-inner automorphism of $R$.   In particular if we let $r^\sigma = y^{-1} ry$, for
all $r \in R$, then we have 
$$
x^\sigma = qx + \alpha z + \beta y + \gamma.
$$  
Since $d(y), d(z) < d(x)$, it is indeed the case that $x^\sigma = qx +$ terms of lower
degree. In Section 4, we will revisit this example and completely determine the group of
$X$-inner automorphisms of $R$.
\vskip.1in

At this point, we turn our attention to enveloping algebras of Lie color algebras and we
begin by defining the appropriate terms.  Let $L$ be a vector space
over the field $K$ and let $G$ be a group.  Then $L$ is a
{\it $G$-graded algebra} if there exist $K$-subspaces $L_g$ such that 
$L = \oplus_{g \in G} L_g$ and a $K$-linear multiplication 
$[\ ,\ ]:L\times L\to L$ satisfying  $[L_g,L_h] \subseteq L_{gh}$ for all $g,h \in G$. 
As usual, the elements of $\bigcup_{g \in G}L_g$ are said to be {\it homogeneous}.
If $G$ is a finite abelian group, then we call a map
$\varepsilon : G \times G \to K^*$ an {\it alternating bicharacter\/} if $\varepsilon (g,hk) =
\varepsilon(g,h) \varepsilon(g,k)$ and
$\varepsilon(g,h) = \varepsilon(h,g)^{-1}$ for all $g,h,k \in G$.

We say that $L$ is a {\it Lie color algebra} over a field $K$ if $L$ is a $G$-graded
algebra and there exists a bicharacter $\varepsilon : G \times G \to K^*$ such that 
$[x,y] = -\varepsilon(g,h) [y,x]$ and $[[x,y],z] = [x,[y,z]] - \varepsilon(g,h) [y,[x,z]]$ 
for all $x \in L_g$, $y \in L_h$, and $z \in L$ (although minor additional technical 
assumptions  are required if the characteristic of $K$ is $2$ or $3$.)

If $g \in G$, then it is easy to check that either $\varepsilon(g,g) = 1$ or $-1$.
Thus, we can partition $G$ into the sets
$G_+ = \{g \in G \mid \varepsilon(g,g) = 1 \}$ and $G_- = \{g \in G \mid
\varepsilon (g,g) = -1 \}$, and we define
$L_+ = \bigcup_{g \in G_+} L_g$ and $L_- = \bigcup_{g \in G_-} L_g$.
If the characteristic of $K$ is $2$, then by convention $G=G_+$.
We say that $L$ is an {\it even Lie color algebra} if $L = L_+$.

Given an even Lie color algebra $L_+$, the Poincar\'e-Birkhoff-Witt Theorem guarantees
that there exists a unique largest $K$-algebra
$U(L_+)$ containing $L_+$ such that $U(L_+)$ is generated by $L_+$ with relations 
$xy - \varepsilon(g,h)yx = [x,y]$, for all $x \in L_g$ and $y \in L_h$.  More precisely, we
have that if $\Cal B$ is a totally ordered basis for $L_+$ consisting of homogeneous elements
then $U(L_+)$ has as a $K$-basis the collection of all ordered monomials
${b_1}^{\beta_1}{b_2}^{\beta_2}\cdots {b_n}^{\beta_n}$, such that $b_i \in {\Cal B}$, 
$b_1 < b_2 < \dots < b_n$ and the $\beta_i$ are nonnegative integers.  
Of course, the grading of $L_+$ by $G$ extends to a grading of $U(L_+)$ and we refer to
homogeneous subsets of $U(L_+)$ in a completely analogous manner.

In light of the Poincar\'e-Birkhoff-Witt Theorem, it is clear that $U(L_+)$ is a domain and if we
filter $U(L_+)$ by letting every homogeneous basis element of $L_+$ have degree $1$, then
$gr(U(L_+))$ is isomorphic to $K_q[X]$.  It is worth noting that if an even Lie color algebra
$L_+$ has trivial multiplication, then $U(L_+)$ is quantum space $K_q[X]$. Conversely, given the
algebra $K_q[X]$ there  always exists a grading group $G$ and bicharacter $\varepsilon$ such that
$K_q[X]$ is the enveloping algebra of an even Lie algebra with trivial multiplication.

Let $T$ be the subgroup of $G$ generated by all the elements of $G$ appearing in the
support of elements of $L_+$.  We can certainly consider $L_+$ and $U(L_+)$ as being
graded by $T$ and is clear that no information is lost when we consider $L_+$ as
a $T$-graded algebra. 
Next, let $H = \{g \in G \mid \varepsilon(g,G) = 1 \}$;  we say that 
$L_+$ is {\it faithfully graded}  by $G$ if $H = 1$. 
The grading of $L_+$ by $G$ certainly induces a grading by the
quotient group ${\overline G} = G/H$.  
Clearly no information about the structure of $L_+$ or $U(L_+)$
is lost when we consider $L_+$ as a ${\overline G}$-graded algebra.

We say that the grading of $L_+$ by $G$ is {\it proper}
if $L_+$ is faithfully graded by $G$ and $G$ is generated by the support of $L_+$.
In light of the comments above, without loss of generality, we may consider even Lie color
algebras where the grading by $G$ is proper.

\proclaim{Theorem 2.5}
Let $L_+$ be an even Lie color algebra with a proper grading by the abelian group $G$.
If $\sigma$ is an $X$-inner automorphism of $U(L_+)$ then there exists some $h \in G$ such that
\roster
\item"(1)"  $x^\sigma = \varepsilon(h,g)x$, for all $x \in L_g$ where $g \not= 1$
\item"(2)"  $x^\sigma = x + \alpha_x$, for all $x \in L_1$ where $\alpha_x$ is a nonzero element
of $K$.
\endroster
\endproclaim

\demo{Proof}
We can filter $U(L_+)$ be letting every homogeneous $x \in L_+$ have degree $1$.  Then
$gr(U(L_+))$ is isomorphic to $K_q[X]$, where $X$ is any homogeneous basis of $L_+$.

The algebra $K_q[X]$ inherits the $G$-grading on $U(L_+)$ and we can extend the grading to $K_q[X,X^{-1}]$.
Let ${\overline \sigma}$ be the $X$-inner automorphism of $K_q[X]$ induced by the action of
$\sigma$ on $U(L_+)$.  Then, by Theorem 2.2,
$$
{\overline \sigma}({\overline x}) = \Delta {\overline x} \Delta^{-1},
$$ 
for all homogeneous $x \in L_+$.

However, $\Delta$ is a homogeneous element of $K_q[X,X^{-1}]$, thus
$\Delta \in K_q[X,X^{-1}]_h$, for some $h \in G$.
Therefore, if $x \in L_g$ we have 
$$
{\overline \sigma}({\overline x}) = \Delta {\overline x} \Delta^{-1} = \varepsilon (h,g) {\overline x} 
\Delta \Delta^{-1} = \varepsilon (h,g) {\overline x}.
$$
In light of this, $\sigma(x) -\varepsilon(h,g)$ has smaller degree than $x$, therefore there
exists some $\beta_x \in K$ such that
$$
\sigma(x) = \varepsilon(h,g)x + \beta_x,
$$
for all $x \in L_g$.

It remains to show that if $g \not=1$ then $\beta_x = 0$.
Since the grading is proper, there is some $k$ is the support of $L_+$ such that
$\varepsilon(g,k) \not=1$.
Let $x \in L_g$, $y \in L_k$ and consider $z = [x,y] = xy - \varepsilon(g,k) yx$.
Applying $\sigma$ to $z$, we obtain
$$\align
z^\sigma&=x^\sigma y^\sigma - \varepsilon(g,k) y^\sigma x^\sigma\\
&=(\varepsilon(h,g)x + \beta_x)(\varepsilon(h,k)y + \beta_y ) -\varepsilon
(g,k)
(\varepsilon(h,k)y + \beta_y )(\varepsilon(h,g)x + \beta_x)\\
&= \varepsilon(h,gk)(xy - \varepsilon(g,k)yx) +\varepsilon(h,g)\beta_y
(1-\varepsilon(g,k))x
+ \varepsilon(h,k)\beta_x (1-\varepsilon(g,k))y\\
&+ \beta_x \beta_y (1-\varepsilon(g,k)).
\endalign
$$

On the other hand, since $z^\sigma = \varepsilon(h,gk)z + \beta_z$, the $k$-component of $z$ is
$0$. Therefore the $k$-component of $x^\sigma y^\sigma - \varepsilon(g,k) y^\sigma x^\sigma$ must
also be $0$ and so, $\varepsilon(h,k)\beta_x (1-\varepsilon(g,k))y = 0$.
However, $\varepsilon (h,k) \not= 0$ and $\varepsilon(g,k) \not= 1$ and it follows that $\beta_x
= 0$, as required.
\qed
\enddemo

We now take a quick look at an example which illustrates the type of $X$-inner automorphism 
described in Theorem 2.5.
\vskip.1in

\noindent {\bf Example 2.6.}
Let $G$ be the free abelian group generated by $g,h$ and define the bicharacter
$\varepsilon$ as $\varepsilon(g,h) = q$, where $q$ is a nonzero element of $K$.
Now let $L_+$ be the $3$-dimensional even Lie color algebra with homogeneous
basis elements $x \in L_1$, $y \in L_g$, $z \in L_h$ such that the only non-trivial
bracket among the basis elements is $[x,y]=y$.
Therefore $U(L_+)$ is the $K$-algebra generated by $x,y,z$ subject to the
relations 
$$
xy - yx = y \ , xz - zx = 0 \ , yz - q zy = 0.
$$
Note that $y$ is a normalizing element of $U(L_+)$ and we can let $\sigma$ be the $X$-inner
automorphism defined as $r^\sigma = y^{-1}ry$, for all $r \in U(L_+)$.
It now follows that
$$
x^\sigma = x+1 \ , y^\sigma = y \ , z^\sigma = q^{-1}z.
$$
In Section 4, we will revisit this example and completely determine its group of
$X$-inner automorphisms.
\vskip.1in

The grading of $G$ on $U(L_+)$ also corresponds to an action of $G$ on $U(L_+)$ as
automorphisms.  Since any automorphism of $U(L_+)$ extends uniquely to $Q = Q(U(L_+))$, this implies
that the grading by $G$ extends to all of $Q$.  Furthermore, every homogeneous 
$x \in L_+$ acts on $U(L_+)$ as a skew derivation and this action also extends uniquely to all
of $Q$.  In light of this, we can now define the {\it homogeneous semi-invariants} in
$Q$ under the adjoint action of $L_+$ as those homogeneous $q \in Q$
such that $\delta_x(q) = \alpha_x q$, for all homogeneous $ x \in L_+$, where 
$\delta_x$ is the skew derivation induced by $x$ and $\alpha_x \in K$.
In particular, this means that if $q \in Q_g$ then for every $x \in L_h$, there is some
$\alpha_x \in K$ such that $xq - \varepsilon(g,h) qx = \alpha_x q$ in $Q$.
Note that if $q$ is a homogeneous semi-invariant then $qU(L_+) = U(L_+)q$, thus $q$ induces an
$X$-inner automorphisms of $U(L_+)$.  We now show that all $X$-inner automorphisms of $U(L_+)$
arise this way.

\proclaim{Theorem 2.7}
{\rm (1)}  The group of $X$-inner automorphisms of $U(L_+)$ is abelian.
\roster
\item"(2)"  The elements in $Q(U(L_+))$ which induce $X$-inner automorphisms of $U(L_+)$ are 
precisely the homogeneous semi-invariants under the adjoint action of $L_+$.
\item"(3)"  The group of homogeneous semi-invariants in $Q(U(L_+))$ is nilpotent of class $2$.
\endroster
\endproclaim

\demo{Proof}
(1)  By Theorem 2.5, it is clear that if $\sigma_1 , \sigma_2$ are $X$-inner automorphisms of 
$U(L_+)$, then $\sigma_1 \sigma_2$ and $\sigma_2 \sigma_1$ agree on all homogeneous elements
of $L_+$.  Since these homogeneous elements generate $U(L_+)$, it follows that 
$\sigma_1 \sigma_2$ and $\sigma_2 \sigma_1$ agree on all of $U(L_+)$, thus the group of
$X$-inner automorphisms of $U(L_+)$ is abelian.

(2)  Let $q \in Q=Q(U(L_+))$ be a homogeneous semi-invariant and let $\sigma$ denote the $X$-inner
automorphism of $U(L_+)$ induced by $q$.  Then, by Theorem 2.5, there  exists some $h \in G$
such that 
$$
qx = \varepsilon(h,g) xq,
\tag 4
$$ 
for all $x \in L_g$ with $g \not=1$ and
$$
qx = xq + \alpha_xq,
\tag 5
$$ 
for all $x \in L_1$.
We claim that $q$ belongs to the $h$-component of $Q$.
To this end, let $q = \sum_{g \in G} q_g$ be the representation of $q$
as a sum of its homogeneous components and suppose $w = q_k$ is nonzero, for some $k \in G$.
Looking at the $kg$-component in equation (4) and the $k$-component in equation (5) 
we see that $wx = \varepsilon(h,g)xw$, for all $x \in L_g$ with $g \not=1$ and
$wx = xw + \alpha_xw$, for all $x \in L_1$. 
As a result, $r^\sigma w = wr$, for all $r \in U(L_+)$, and therefore $\sigma$
is also induced by $w$.

If $x \in L_g$ with $g \not=1$ and if $\delta$ denotes the skew derivation induced by $x$, then
there exists a nonzero $a \in U(L_+)$ such that $wa, \delta(w)a$ both belong to $U(L_+)$.
Since 
$$\delta(wa) = \delta(w)a + \varepsilon(g,k)w\delta(a),
\tag 6
$$ it follows that
$w\delta(a) \in U(L_+)$.  
If we let $d$ denote the degree of an element in $U(L_+)$, then we can use the fact that
$(wa)\delta(a) = (a^\sigma w)(\delta(a)) =  a^\sigma(w\delta(a))$, to conclude that
$d(wa) + d(\delta(a)) = d(a^\sigma) + d(w\delta(a))$.
However, we know that $\sigma$ preserves degree, thus
$d(a^\sigma) = d(a) \ge d(\delta(a))$.
Therefore $d(wa) \ge d(w\delta(a))$.

If we let $M$ denote $d(wa)$, then we see that $d(\delta(wa)),d(w\delta(a)) \le M$.
It now follows using equation (6) that $d(\delta(w)a) \le M$.
Since $w$ induces $\sigma$, we have $wx = \varepsilon(h,g)xw$.
On the other hand, $\delta(w) = xw -\varepsilon(g,k)wx$ and if we combine the previous two
equations we see that  
$$
\delta(w) = xw - \varepsilon(g,k)\varepsilon(h,g)xw = (1-\varepsilon(g,kh^{-1}))xw.
$$
If $h \not=k$, we could have chosen $g$ such that $\varepsilon(g,kh^{-1}) \not=1$.
However, if we multiply the previous equation on the right by $a$ this now implies
that $M \ge d(\delta(w)a) = d(xwa) = d(x) + d(wa) = M+1$, a contradiction.
Therefore $h=k$ and it follows that $q$ is indeed homogeneous as it belongs to the $h$-component
of $Q$.

Finally, if $x \in L_g$ with $g \not=1$, then we have
$$
\delta(q) = xq - \varepsilon(g,h)qx = \varepsilon(g,h)(\varepsilon(h,g)xq - qx) = 
\varepsilon(g,h)(x^\sigma q - qx) = 0
$$ 
and if $x \in L_1$
then $\delta(q) = xq - qx = -\alpha_x q$.
Thus $q$ is a homogeneous semi-invariant.

(3)  If $q_1, q_2 \in Q$ are homogeneous semi-invariants then, by (1), they induce
commuting $X$-inner automorphisms.  Thus $q_2^{-1}q_1^{-1}q_2q_1$ is in the center
of $Q$ and it follows that the group of homogeneous semi-invariants is nilpotent of
class $2$.
\qed
\enddemo

We conclude this section by using Theorem 2.5 to study crossed products.

\proclaim{Theorem 2.8}
Let $L_+$ be an even Lie color algebra over a field of characteristic $0$
such that the every value of the bicharacter $\varepsilon$ is either $1$ or not a root of
$1$.  If $G$ is a finite group of automorphisms of $U(L_+)$, then any crossed product
$U(L_+)*_tG$ is prime.  Furthermore, if $U(L_+)$ is primitive then so is $U(L_+)*_tG$.
\endproclaim

\demo{Proof}
Let $\sigma$ be an $X$-inner automorphism of $U(L_+)$ and let $x$ be a homogeneous element of $L_+$.
It follows from Theorem 2.5 that if $\sigma(x) \not= x$, then $\sigma^n(x) \not= x$, for any 
positive integer $n$.  Thus any non-identity $X$-inner automorphism of $U(L_+)$ must have
infinite order.

If $G$ is a finite group of automorphisms of $U(L_+)$, then no element of $G$ other than the
identity is $X$-inner.  As a result, every nonzero ideal of a crossed product $U(L_+)*_tG$ must
intersect $U(L_+)$ nontrivially.  Since $U(L_+)$ is a domain, $U(L_+)*_tG$ must be prime.
In addition, since $U(L_+)*_tG$ is a free module of finite rank over $U(L_+)$, a result in
\cite{B} implies that if $U(L_+)$ is primitive then so is $U(L_+)*_tG$.
\qed
\enddemo

\head
\S 3. $q$-Skew Polynomial Rings
\endhead

In \cite{LM}, Leroy and Matczuk examine the $R$-stabilizing $X$-inner
automorphisms of skew polynomial rings $R[x;\tau,\delta]$, where $R$ is a prime ring. 
Recall that an $X$-inner automorphism $\sigma$ of $R[x;\tau,\delta]$ is $R$-stabilizing if
$R^\sigma = R$.  It is not hard to show that when $R$ is a domain, every $X$-inner
automorphism of $R[x;\tau,\delta]$ is $R$-stabilizing. In \cite{LM} we are introduced to
{\it monic invariant polynomials}, which are monic polynomials $p(x)$ in
$Q(R)[x;\tau,\delta]$ with the properties that (1) $p(x)x = (ax+b)p(x)$, for some $a,b \in
Q(R)$ and (2) for every $r \in Q(R)$ there is an $s \in Q(R)$ such that  $p(x)r = sp(x)$.  
We will let $P(x)$ denote the unique monic invariant polynomial in $Q(R)[x;\tau,\delta]$ of
minimal nonzero degree if such a polynomial exists and let $P(x)=1$ if no such polynomial
exists.

Theorem 4.3 in \cite{LM} states that every $R$-stabilizing $X$-inner automorphism of 
$R[x;\tau,\delta]$ is induced by an element of the form $P(x)^mw$, where $m$ is an integer
and $w \in N(R)$.  Leroy and Matczuk show that if $\delta$ and $\tau$ commute, then $w$ and
$P(x)$ both induce $X$-inner automorphisms of $R[x;\tau,\delta]$.  
However, in general, this is not the case as we shall see in Section 4 when we study the
quantum Weyl algebra. In this section, we will refine the result in \cite{LM} in the special
case where $\delta$ is a $q$-skew derivation.  This will help us to determine 
the group of $X$-inner automorphisms of various semi-commutative quantum algebras in
Section 4.

We say that $R[x;\tau,\delta]$ has the {\it $(R,R)$-bimodule property} if every nonzero
$(R,R)$-bimodule of $R$ intersects $R$ nontrivially.

\proclaim{Proposition 3.1}
Let $R$ be a prime ring with a $q$-skew derivation $\delta$ such that either
$q$ is not a root of $1$ or $q=1$ where $R$ has characteristic $0$.
Then either  $R[x;\tau,\delta]$ has the $(R,R)$-bimodule property or $\delta$ is $X$-inner.
\endproclaim

\demo{Proof}
Suppose $I$ is an $(R,R)$-bimodule of $R[x;\tau,\delta]$ which does not intersect $R$
nontrivially; we will show that $\delta$ must be $X$-inner.  
Let $w = c_nx^n + \cdots + c_1x + c_0$ be a nonzero element of $I$ of minimal degree $n$.
Next, let 
$$
J = \{ r \in R \mid I \ \text{contains an element of degree} \ n \ 
\text{with leading coefficient} \ r\}.
$$
Since $w$ has minimal degree, for every $r \in J$, there is a unique element
$u \in I$ of the form 
$$
u = rx^n + r_{n-1}x^{n-1} + \cdots + r_1x + r_0.
$$
We can now define the left-module maps $b_i: J \longrightarrow R$, for $0 \le i \le n-1$,
as $r \cdot b_i = r_i$.  The maps $b_i$ belong to the left Martindale quotient ring $Q_l$ of
$R$ and both $\tau$ and $\delta$ extend uniquely to $Q_l$.
Let 
$$
v = x_n + b_{n-1}x^{n-1} + \cdots + b_1x + b_0 \in Q_l[x;\tau,\delta];
$$
it is clear that $Jv \subseteq I$.
If $s \in R$, consider 
$$
z = \tau^n(s)v - vs = (\tau^n(s) - \tau^n(s))x^n + \text{terms of lower degree}.  
$$
Note that the degree of $z$ is less than $n$ and if $z \neq 0$  then 
$Jz = J(\tau^n(s)v - vs)$ would contain elements of $I$ of degree less than $n$.
Thus $z=0$.

If we examine the coefficient of the degree $n-1$ term of $z$ we see that
$$
0 = \tau^n(s) b_{n-1} - b_{n-1}\tau^{n-1}(s) - (1 + q^{-1} + \cdots +
q^{n-1})\delta(\tau^{n-1}(s)).
$$
Letting $\alpha = -( 1 + q^{-1} + \cdots + q^{n-1})$, $t = \tau^{n-1}(s)$, and $b = b_{n-1}$ 
in the previous equation, we obtain
$$
\alpha\delta(t) = bt - \tau(t) b,
$$
for all $t \in R$.
Note that since either $q$ is not a root of $1$ or $q=1$ where $R$ has characteristic $0$,
it follows that $\alpha$ is a nonzero element of $K$.
Thus $\delta$ is $X$-inner as it is induced by the element $\alpha^{-1}b$ which also 
lives in $Q(R)$.
\qed
\enddemo

We can now combine Proposition 3.1 with the result from \cite{LM} to describe the
$R$-stabilizing $X$-inner automorphisms when $\delta$ is a $q$-skew derivation.
Before proving Theorem 3.2, we should point out that the assumption that $\delta$ 
is $q$-skew is only used to obtain the dichotomy that either $R[x;\tau,\delta]$
has the $(R,R)$-bimodule property or $\delta$ is $X$-inner.
If $\delta$ is not $q$-skew but one already knows that 
either $R[x;\tau,\delta]$ has the $(R,R)$-bimodule property or $\delta$ is
$X$-inner, then one can still use the appropriate part of the theorem below.

\proclaim{Theorem 3.2}
Let $R$ be a prime ring with a $q$-skew derivation $\delta$ such that either $q$ is not
a root of $1$ or $q=1$ where $R$ has characteristic $0$.
If $\sigma$ is an $R$-stabilizing $X$-inner automorphism of $R[x;\tau,\delta]$, then either
\roster
\item"(1)"  $\sigma$ is induced by some $w \in N(R)$.   Furthermore, an element $w \in N(R)$
induces an $X$-inner automorphism of $R[x;\tau,\delta]$ if and only if
$w^{-1}\tau(w)$, $\tau(w^{-1})w$, and $w^{-1}\delta(w)$  belong to $R$.  In this case,
$$
x^\sigma = w^{-1}\tau(w)x + w^{-1}\delta(w).
$$  
This case occurs when $R[x;\tau,\delta]$ has the $(R,R)$-bimodule property.
\endroster
or
\roster
\item"(2)" $\sigma$ is induced by an element of the form $(x-c)^mw$, where
$w \in N(R)$, $m$ is an integer, and $\delta$ is $X$-inner and induced by $c$.
Furthermore, an element of the form $(x-c)^mw$, with $w,m,c$ as above, induces an
$X$-inner automorphism of $R[x;\tau,\delta]$ if and only if
$w^{-1}\tau(w), \tau(w^{-1})w, \tau^{-m}(c) - \tau(w)cw^{-1} \in R$.  In this case,
$$
r^\sigma = w^{-1}\tau^m(r)w,
$$ 
for all $r \in R$, and
$$
x^\sigma = w^{-1}\tau(w)x + w^{-1}(\tau^{-m}(c) - \tau(w)cw^{-1})w.
$$  
This case occurs when $\delta$ is $X$-inner.
\endroster
\endproclaim

\demo{Proof}
If $p(x)$ is a monic invariant polynomial, then 
$Q(R)p(x) \cap R[x;\tau,\delta]$ is  a nonzero $(R,R)$-bimodule
which does not intersect $R$ nontrivially.  Therefore, if $R[x;\tau,\delta]$ has the
$(R,R)$-bimodule property it must be the case that there is no monic irreducible
polynomial of minimal degree and so, $P(x)=1$.  However, in this case Theorem 4.3 of
\cite{LM} implies that there exists some $w \in N(R)$ such that
$$
s^\sigma = w^{-1}sw,$$
for all $s \in R[x;\tau,\delta]$, where
$w^{-1}\tau(w), w^{-1}\delta(w) \in R$ and 
$$
x^\sigma = w^{-1}\tau(w)x + w^{-1}\delta(w).
$$
Furthermore, since $\sigma$ is surjective it must also be the case that 
$(w^{-1}\tau(w))^{-1} = \tau(w^{-1})w \in R$.

On the other hand, if $\delta$ is $X$-inner then there exists some $c \in Q(R)$ such that
$\delta(r) = cr - \tau(r)c$, for all $r \in Q(R)$.  However, since we already know that
$\delta(r) = xr - \tau(r)x$, for all $r \in Q(R)$, it follows that
$(x-c)r = \tau(r)(x-c)$, for all $r \in Q(R)$.  
In addition, we have
$$
(x-c)x = (x-c)x + (-\delta(c)) - (-\delta (c)) = 
$$
$$
(x^2-cx) + (\tau(c)x - cx) - 
(\tau(c)c - c^2)  = (x + \tau(c) - c)(x-c).
$$
Thus $x-c$ is a monic irreducible polynomial of minimal degree and, by Proposition 3.4 of
\cite{LM}, it is the unique one of that degree.   
Hence $P(x) = x-c$. 
When we let $P(x) = x-c$ in Theorem 4.3 of \cite{LM}, its conclusion  can
be restated to say that there exists some $w \in N(R)$ and integer $m$ such that  
$$
s^\sigma = (P(x)^mw)^{-1}s(P(x)^mw),
$$ 
for all $s \in R[x;\tau,\delta]$, where 
$w^{-1}\tau(w), \tau^{-m}(c) - \tau(w)cw^{-1} \in R$.
Furthermore, 
$$
r^\sigma = w^{-1}\tau^m(r)w,
$$ 
for all $r \in R$, and 
$$
x^\sigma = w^{-1}\tau(w)x + w^{-1}(\tau^{-m}(c) + \tau(w)cw^{-1})w.
$$
As in the previous case, it must also be true that  
$(w^{-1}\tau(w))^{-1} = \tau(w^{-1})w \in R$.
 
It remains to determine which elements of the form $w$ in part (1) and
$P(x)^mw$  in part (2) do indeed induce $X$-inner automorphisms of $R[x;\tau,\delta]$.
By Lemma 4.2 of \cite{LM}, the $X$-inner automorphisms of $Q(R)[x;\tau,\delta]$ induced
by $w$ and $P(x)^mw$ will be $X$-inner automorphisms of $R[x;\tau,\delta]$ provided they
stabilize $R[x;\tau,\delta]$.  
If we let $v=w$ or $v=P(x)^mw$, depending upon whether we are in part (1) or part (2), then
we know that $Rv = vR$.  Therefore in order to complete the proof and show that
$R[x;\tau,\delta]v = vR[x;\tau,\delta]$, it suffices to show that 
$xv = v(ax+b)$, where $a, a^{-1},b \in R$.

For part (1), note that since $\delta(w) = xw - \tau(w)x$, it follows that
$$
xw = \tau(w)x + \delta(w) = w(w^{-1}\tau(w) + w^{-1}\delta(w)).
$$
Since $w^{-1}\tau(w), \tau(w^{-1})w, w^{-1}\delta(w) \in R$,
it is the case that $xw = w(ax+b)$, where $a, a^{-1},b \in R$.
This complete the proof of part (1).

To complete part (2), we consider the case where $v=P(x)^mw$.
We now have
$$
(x-c)v = (x-c)^{m+1}w = (x-c)^m\tau(w)(x-c) = 
$$
$$
((x-c)^mw)(w^{-1}\tau(w))(x-c) =
v(w^{-1}\tau(w))(x-c).
\tag 7
$$
Furthermore, we also have
$$
cv = c(x-c)^mw = (x-c)^m\tau^{-m}(c)w = 
$$
$$
((x-c)^mw)(w^{-1}\tau^{-m}(c)w) 
= v(w^{-1}\tau^{-m}(c)w).
\tag 8
$$
If we add equations (7) and (8), we obtain
$$
xv = v(w^{-1}\tau(w)(x-c) + w^{-1}\tau^{-m}(c)w) = v(w^{-1}\tau(w)x + w^{-1}\tau^{-m}(c)w
- w^{-1}\tau(w)c).
$$
However 
$$
w^{-1}\tau^{-m}(c)w - w^{-1}\tau(w)c = w^{-1}(\tau^{-m}(c) - \tau(w)cw^{-1})w
$$
and this belongs to $R$ if and only if $\tau^{-m}(c) - \tau(w)cw^{-1} \in R$.
Therefore, it is once again the case that
$xv = v(ax+b)$, where $a, a^{-1},b \in R$, thereby completing the proof of part (2).
\qed
\enddemo

Note that although $P(x)$ and $w$ induce $X$-inner automorphisms of $Q(R)[x;\tau,\delta]$,
they do not necessarily induce $X$-inner automorphisms of $R[x;\tau,\delta]$ even if
$P(x)^mw$ induces an $X$-inner automorphism of $R[x;\tau,\delta]$.
We shall see an example of this in the next section when we compute
the $X$-inner automorphisms of the quantum Weyl algebra.

\head
\S 4.  Examples
\endhead

In this section we will apply the results of the previous two sections to determine the
group of $X$-inner automorphisms of various semi-commutative quantum algebras.
\vskip.1in

\noindent {\bf Example 4.1.}
Let $R$ be the quantum Weyl algebra $K[x,y \mid xy - qyx =1]$, where $q$ is an element of
the field $K$ which is not a root of $1$.  If $\sigma$ is an $X$-inner automorphism of
$R$ and if we filter $R$ by letting $x$ and $y$ both have degree $1$, then Theorem 2.3
implies that 
$$
x^\sigma = \alpha_1x + \beta_1 \ \text{and} \  y^\sigma = \alpha_2y + \beta_2,
$$ 
where $\alpha_1, \alpha_2$ are both powers of $q$ and $\beta_1, \beta_2 \in K$.
However, since $x^\sigma y^\sigma - q y^\sigma x^\sigma = 1$, it follows that
$\beta_1 = \beta_2 = 0$ and $\alpha_2 = \alpha_1^{-1}$.
Therefore, we now have
$$
x^\sigma = q^n x \  \text{and} \  y^\sigma = q^{-n}y,
$$
where $n$ is an integer.

We can also view $R$ as the skew polynomial ring $K[y][x;\tau, \delta]$ where
$\tau(y) = qy$ and   $\delta(y) = 1$.  Note that $\delta$ is $X$-inner as it
is induced by the element $c = (1-q)^{-1}y^{-1}$.
Since $\delta$ is $X$-inner, we can apply part (2) of Theorem 3.2 to conclude that
$\sigma$ is induced by an element of the form
$P(x)^mw$, where $P(x) = x-c$, $m$ is an integer, and $w \in N(K[y])$.
Furthermore, it must also be the case that
$$
w^{-1}\tau(w) = q^n \  \text{and} \  \tau^{-m}(c) - \tau(w)cw^{-1} = 0.
$$
The first equation from the previous line implies that $w$ must be a scalar 
multiple of $y^n$ therefore, without loss, we may assume $w = y^n$.
On the other hand, the second equation from above implies that 
$q^mc - q^nc = 0$, and so $m=n$.
Therefore $\sigma$ must be induced by an element of the form $P(x)^ny^n$.

Observe that 
$$P(x)y = (x-(1-q)^{-1}y^{-1})y = (xy -(1-q)^{-1}) =  
(1-q)^{-1}((1-q)xy -1) = 
$$
$$
(1-q)^{-1}(-qxy + xy -1) = (1-q)^{-1}(-qxy +(qyx+1) -1)
=q(1-q)^{-1}(xy-yx).
$$
In light of this, $xy-yx$ induces an $X$-inner automorphism of $R$ and it  follows that
$P(x)^ny^n = \gamma (xy-yx)^n$, where $\gamma$ is a nonzero element of $K$.
Therefore all the $X$-inner automorphisms of $R$ are induced by powers of $xy-yx$ and the
$X$-inner automorphism group of $R$ is infinite cyclic.
\vskip.1in

In the above example, note that although $P(x)^ny^n$ induces an $X$-inner automorphism of
$R$, for every $n$, neither $y$ nor $P(x)$ induces an $X$-inner automorphism of $R$.
In particular, we have
$$
y^{-1}xy = qx+y^{-1} \notin R \  \text{and} \  P(x)^{-1}yP(x) = x - y^{-1} \notin R.
$$
In Theorem 4.7 of \cite{LM}, it shown that if $\delta$ and $\tau$ commute then the elements
$P(x)$ and $w$ described in Theorem 3.2 both induce $X$-inner automorphisms of $R$.
However in the quantum Weyl algebra example above, this is not the case as 
$\delta \tau  = q \tau \delta$.

\vskip.1in
\noindent {\bf Example 4..2. (or Example 2.6 revisited)}
Let $G$ be the free abelian group generated by $g,h$ with bicharacter
$\varepsilon$ defined as $\varepsilon(g,h) = q$, where we now assume that $q \in K$ is 
not a root of $1$. 
Let $L_+$ be the $3$-dimensional even Lie color algebra with homogeneous
basis elements $x \in L_1$, $y \in L_g$, $z \in L_h$ such that the only non-trivial
bracket among the basis elements is $[x,y]=y$.
Therefore $U(L_+)$ is the $K$-algebra generated by $x,y,z$ subject to the
relations 
$$
xy - yx = y,   \  xz - zx = 0,   \  yz - q zy = 0.
$$

Observe that $y$ and $z$ are homogeneous semi-invariants of $U(L_+)$. 
Therefore, by Theorem 2.7, it follows that every monomial of the form $y^nz^m$ induces an
$X$-inner automorphism of $U(L_+)$, where $n,m$ are integers.
We need to show that these are the only $X$-inner automorphisms of $U(L_+)$. 

We can consider $U(L_+)$ to be the differential operator ring  $K_q[y,z][x;d]$, where
$K_q[y,z]$ is the quantum plane and $d(y) = y$ and $d(z) = 0$.
Since $q$ is not a root of $1$, the extended centroid of $K_q[y,z]$ is the field $K$ and 
Lemma 1.1(b) now implies that every nonzero ideal of $K_q[y,z]$ contains a monomial of the form
$y^iz^j$.  Therefore the Martindale quotient ring $Q(K_q[y,z])$ is obtained from
$K_q[y,z]$ by inverting monomials. 

If the derivation $d$ is $X$-inner, then it is induced by some 
$c = \sum \alpha_{i,j} y^iz^j$, where $i,j$ are integers.  
Since $d(z) = 0$, $c$ and $z$ commute
and it follows that $c$ must be of the form $\sum \beta_j z^j$.
However, we also know that $d(y) = y$, which implies that
$$
y = d(y) = cy-yc = \sum \beta_j (q^{-j}-1)yz^j,
$$
which is impossible.
Hence $d$ is not $X$-inner and we can apply part (1) of Theorem 3.2.

Therefore, if $\sigma$ is an $X$-inner derivation of $U(L_+)$, then $\sigma$ must be generated
by a normalizing element in $Q(K_q[y,z])$.  However, Proposition 2.2 implies that the only
normalizing elements are  monomials, thus $\sigma$ is indeed induced by $y^nz^m$, for some
integers $n,m$. 

If we define the $X$-inner automorphism 
$\sigma_{n,m}$ as $\sigma_{n,m}(r) = (y^nz^m)^{-1}r(y^nz^m)$, for all $r \in U(L_+)$, we have 
$$
\sigma_{n,m}(x) = x+n, \    \sigma_{n,m}(y) =  q^my, \     \sigma_{n,m}(z) = q^{-n}z.
$$
Since every pair of integers $(n,m)$ induces a different automorphism, the
group of $X$-inner automorphisms of $U(L_+)$ is free abelian on
two generators. 
\vskip.1in

We now examine the first of two $q$-analogs of the Heisenberg algebra studied in \cite{KS}.

\vskip.1in
\noindent {\bf Example 4.3.}  Let $R$ be the $K$-algebra generated by $x,y,z$ 
subject to the
relations
$$
xy-qyx = z, \   xz-qxz = 0, \  yz - q^{-1}zy = 0,
$$
where $q$ is an element of the field $K$ which is not a root of $1$
and $K$ has characteristic $0$.
We can view $R$ as $R = K_q[y,z][x;\tau, \delta]$, where
$ K_q[y,z]$ is the quantum plane and
$$
\delta(y) = z, \  \delta(z) = 0, \  \tau(y) = qy, \  \tau(z) = qz.
$$

By Theorem 2.3, 
$$
x^\sigma = \alpha_1 x + \beta_1, \  y^\sigma = \alpha_2 y + \beta_2, \  
z^\sigma = \alpha_3 z + \beta_3,
$$
where each $\alpha_i$ is a power of $q$ and each $\beta_i \in K$.
Since $x^\sigma y^\sigma - q y^\sigma x^\sigma = z^\sigma$, it follows that
$\beta_1 = 0$.
As a result, $x^\sigma = q^nx$, for some integer $n$.

Since $q$ is not a root of $1$, $Q(K_q[y,z])$ is obtained
from  $K_q[y,z]$ by inverting monomials.
Therefore, if $\delta$ is $X$-inner then it is induced by an element 
$c = \sum \alpha_{i,j} y^iz^j$, where $i,j$ are integers.
However, we then have
$$
0 = \delta(z) = cz - qzc = \sum \alpha_{i,j}(1-q^{i+1})y^iz^{j+1},
$$
which implies that $c$ must be of the form $\sum \beta_jy^{-1}z^j$.
We next observe that
$$
z = \delta(y) = cy -qyc = \sum \beta_j(q^j-q)z^j,
$$
which is impossible, thus $\delta$ is not $X$-inner.

Since $\delta$ is not $X$-inner and $\delta \tau = \tau \delta$, 
Proposition 3.1 implies that every nonzero $(K_q[y,z],K_q[y,z])$-bimodule of 
$R = K_q[y,z][x;\tau, \delta]$ intersects $K_q[y,z]$ nontrivially.
Therefore, we can apply part (1) of Theorem 3.2 to conclude that $\sigma$
is induced by an element of $Q(K_q[y,z])$ which normalizes $K_q[y,z]$.
However, the extended center of $K_q[y,z]$ is the field $K$ and
Proposition 2.2 implies that the only elements in $Q(K_q[y,z])$ which normalize
$K_q[y,z]$ are monomials.
Thus part (1) of Theorem 3.2 tells us that $\sigma$ is induced by some monomial $w = y^iz^j$, 
where $i,j$ are integers and $w^{-1}\tau(w) = q^n$ and $w^{-1}\delta(w) = 0$.
Thus $n=i$ and we observe that
$$
0 = \delta(w) = \delta(y^nz^j) = \delta(y^n)z^j = nq^{n-1}y^{n-1}z^{j+1}.
$$

Since $K$ has characteristic $0$,  $n=0$ and $\sigma$ must be induced by a power of $z$.
If we let $\sigma_m$ be defined as $\sigma_m(r) = z^{-m}rz^m$, for all $r \in R$, then we 
have
$$  
\sigma_m(x) = q^mx, \     \sigma_m(y) = q^{-m}y , \    \sigma_m(z) = z.
$$
Since each integer $m$ yields a different automorphism, we see that in this case the group of
$X$-inner automorphisms of $R$ is infinite cyclic.
\vskip.1in

We now consider the second $q$-analog of the Heisenberg algebra studied in \cite{KS}.  It is
worth noting that although the generators and relations of these two algebras are almost
identical, the group of $X$-inner automorphisms of the first algebra is infinite
cyclic, whereas the group of $X$-inner automorphisms of
the second algebra is free abelian on two generators.

\vskip.1in
\noindent {\bf Example 4.4.}
Let $R$ be the $K$-algebra generated by $x,y,z$ subject to the relations
$$
xy-qyx = z, \   xz-q^{-1}xz = 0, \  yz - qzy = 0,
$$
where $q$ is an element of the field $K$ which is not a root of $1$
We can view $R$ as $R = K_q[y,z][x;\tau, \delta]$, where
$ K_q[y,z]$ is the quantum plane and
$$
\delta(y) = z, \  \delta(z) = 0, \  \tau(y) = qy, \  \tau(z) = q^{-1}z.
$$

A minor difference between this algebra and the one in the previous example is
that $\delta$ and $\tau$ no longer commute as we have $\delta \tau = q^2 \tau \delta$.
However, a more significant difference is that $\delta$ is now $X$-inner as it is
induced by $c = q(1-q^2)^{-1}y^{-1}z$.

As in the previous example, we can use Theorem 2.3 and the fact that
$x^\sigma y^\sigma - q y^\sigma x^\sigma = z^\sigma$ to conclude that 
$x^\sigma = q^nx$, for some integer $n$.
It is again the case that the only elements in $Q(K_q[y,z])$ which normalize
$K_q[y,z]$ are monomials.  Combining this with the fact that $\delta$ is $X$-inner, 
we can use part (2) of Theorem 3.2 to see that $\sigma$ is induced by an
element of the form $P(x)^mw$, where $P(x)=x-c$ and $w$ must be of the form $y^iz^j$, for 
integers $i,j$.
Furthermore, Theorem 3.2 also tells us  that
$w^{-1}\tau(w) = q^n$ and $\tau^{-m}(c) - \tau(w)cw^{-1} = 0$.

Using the facts that $c = q(1-q^2)^{-1}y^{-1}z$, $\tau(c) = q^{-2}c$ and $w = y^iz^j$, the
second equation from the previous line implies that 
$$
0 = q^{2m}y^{-1}z - \tau(y^iz^j)y^{-1}z(y^iz^j)^{-1} = q^{2m}y^{-1}z - q^{2i}y^{-1}z.
$$
As a result, $m=i$ and every $X$-inner automorphism of $R$ is of the form $\sigma_{n,m}$, where
$\sigma_{n,m}(r) = (P(x)^my^mz^n)^{-1}r(P(x)^my^mz^n)$, for all $r \in R$ and integers $n,m$.
More precisely, we have 
$$
\sigma_{n,m}(x) = q^{m-n}x, \    \sigma_{n,m}(y) = q^{n+m}y  , \    \sigma_{n,m}(z) = q^{-2m}z.  
$$
Since all of the $\sigma_{n,m}$ commute and every pair of integers $n,m$ results in a different
automorphism, we see that the group of $X$-inner automorphisms of $R$ is
free abelian on two generators.
\vskip.1in

Note that in the previous example, despite the fact that $P(x)^my^mz^n$ always induces an $X$-inner
automorphism of $R$, $P(x)$ and $y$ do not induce $X$-inner automorphisms of $R$.
In particular, we have
$$
P(x)^{-1}xP(x) = x - qy^{-1}z  \notin R \   \text{and}  \    y^{-1}xy   = qx + y^{-1}z\notin  R.
$$

We now revisit an algebra which does not resemble a $q$-analog of a Weyl
algebra or enveloping algebra but is semi-commutative and therefore can be studied using the techniques developed in the previous two sections.

\vskip.1in
\noindent {\bf Example 4.5. (or Example 2.4 revisited)}
Let $R = K[x,y,z \mid xy-qyx = \alpha yz + \beta y^2 + \gamma y]$, where $q, \alpha, \beta,
\gamma$ belong to the field $K$ and $q$ is not a root of $1$.  As we saw in Example 2.4, if we
filter $R$ by letting $x$ have degree $2$ and letting $y$ and $z$ have degree $1$, then
$gr(R) = K[\overline x, \overline y, \overline z \mid \overline x \overline y = q \overline y 
\overline x]$.
Therefore, by Theorem 2.3, if $\sigma$ is any $X$-inner automorphism of $R$ then
$x^\sigma = q^n +$ terms of lower degree, for some integer $n$.

Note that $R = K[y,z][x;\tau, \delta]$, where
$$
\delta(y) = \alpha yz + \beta y^2 + \gamma y, \  \delta(z) = 0, \  
\tau(y) = qy, \  \tau(z) = z.
$$
Also note that $K[y,z]$ is commutative and $\delta$ is $X$-inner as it is induced by
$c = (1-q)^{-1}(\alpha z + \beta y + \gamma)$.
Therefore, by part (2) of Theorem 3.2, if $P(x)=x-c$ then $\sigma$ is induced by an element
of the form $P(x)^mw$, where $w^{-1}\tau(x) = q^n$.
However this implies that $w = f(z)y^n$, for some rational function $f(z) \in K(z)$.
Since $z$ is central, without loss of generality, we may assume that $w=y^n$.

Observe that
$$
y^{-1}xy = qx + \alpha z + \beta y + \gamma \in R \  \text{and} \  
P(x)^{-1}xP(x) = x + q^{-1}\beta y \in R.
$$
Therefore both $y$ and $P(x)$ induce $X$-inner automorphisms of $R$ and it is easy to check
these automorphisms commute.
As a result, for every pair of integers $(n,m)$, conjugation by $P(x)^my^n$ does indeed induce
an $X$-inner automorphisms of $R$.
If we let $\sigma_{n,m}$ be defined as
$\sigma_{n,m}(r) = (P(x)^my^n)^{-1}r(P(x)^my^n)$, for all $r \in R$, then we have
$$
\sigma_{n,m}(x) = q^nx + (q^n-1)(q-1)^{-1}\alpha z + (q^n - q^{-m})(q-1)^{-1}\beta y +  
(q^n-1)(q-1)^{-1}\gamma,
$$
$$
\sigma_{n,m}(y) = q^{-m}y, \      \sigma_{n,m}(z) = z.
$$
Since every pair of integers $(n,m)$ induces a different automorphism, the
group of $X$-inner automorphisms of $R$ is free abelian on
two generators.
\vskip.1in

Note that in all of the examples in this section, every non-identity $X$-inner automorphism has infinite order.  In light of this, using the identical argument used in the proof of Theorem 2.8, we have

\proclaim{Corollary 4.6}
Let $R$ be any of the algebras in Examples 4.1-4.5 and let $G$ be a finite group of automorphisms of $R$.  Then any crossed product $R*_tG$ must be prime.
\endproclaim


\Refs
\widestnumber\key{OsbP}

\ref\key B\by J. Bergen
\paper A note on the primitivity of ring extensions
\jour Comm. Algebra 	\vol	23 \yr 1995 
\pages 4625--4631 \endref


\ref\key KS\by E. E. Kirkman and L. W. Small
\paper $q$-analogs of harmonic oscillators and related rings
\jour Israel J. Math. \vol 81 \yr 1993
\pages 111--127 \endref

\ref\key LM\by A. Leroy and J. Matczuk
\paper The extended centroid and $X$-inner automorphisms of Ore extensions
\jour J. Algebra \vol 145 \yr 1992 
\pages 143--177 \endref

\ref\key Mc\by J. C. McConnell 
\paper Quantum groups, filtered rings and Gelfand-Kirillov dimension
\paperinfo Noncommutative Ring Theory (Athens, Ohio, 1989), 
Lecture Notes in Mathematics\vol 1448
\publ Springer-Verlag\publaddr Berlin \yr 1990
\pages 139--147 \endref


\ref\key McP\by J. C. McConnell and J. J. Pettit
\paper Crossed products and  multiplicative analogs of Weyl algebras
\jour J. London Math. Soc. \vol 38 \yr 1988
\pages 47--55 \endref


\ref\key M\by S. Montgomery
\paper $X$-inner automorphisms of filtered algebras
\jour Proc. Amer. Math. Soc.  \vol 83 \yr 1981
\pages 263--268
\endref


\ref\key Mu\by I. Musson
\paper Ring theoretic properties of the coordinate rings of quantum
symplectic and Euclidean space
\paperinfo Ring Theory (Proceedings of the 1992 Biennial Ohio State-Denison Conference)
\publ World Scientific\publaddr Singapore \yr 1993
\pages 248--258 \endref

\ref\key OsbP\by J. M. Osborn and D. S. Passman
\paper  Derivations of skew polynomial rings
\jour J. Algebra \vol 176   \yr 1995
\pages  417--448
\endref


\ref\key OP\by J. Osterburg and D. S. Passman
\paper  $X$-inner automorphisms of enveloping rings
\jour  J. Algebra \vol 130   \yr 1990
\pages 412--434
\endref

\endRefs
\enddocument